\documentclass[12pt,a4paper]{article}
\usepackage{amsmath,amssymb,fancyhdr,url}
\usepackage{graphicx}
\usepackage{multicol}
\usepackage{multirow}
\usepackage{algorithmic}
\usepackage{algorithm}
\usepackage{listings}
\title{An alternative derivation of a new Lanczos-type algorithm for systems of linear equations}
\author{Saifullah,\footnote{Department of Mathematics, University of Peshawar, Khyber Pakhtunkhwa, 25120, Pakistan. E-mail: saifullah.maths@gmail.com} \mbox{ } Muhammad Farooq\footnote{Department of Mathematics, University of Peshawar, Khyber Pakhtunkhwa, 25120, Pakistan. E-mail: mfarooq@upesh.edu.pk} \mbox{ } and Abdellah Salhi\footnote{Department of Mathematical Sciences, University of Essex, Wivenhoe Park, Colchester, CO4 3SQ, UK. E-mail: as@essex.ac.uk}}
\begin{document}
\maketitle
\begin{abstract}
Various recurrence relations between formal orthogonal polynomials can be used to derive Lanczos-type algorithms. In this paper, we consider  recurrence relation $A_{12}$ for the choice  $U_i(x)=P_i(x)$, where $U_i$ is an auxiliary family of polynomials of exact degree $i$. It leads to a Lanczos-type algorithm that shows superior stability when compared to existing Lanczos-type algorithms. The new algorithm is derived and described. It is then computationally compared to the most robust algorithms of this type, namely $A_{12}$, $A_5/B_{10}$ and $A_8/B_{10}$, on the same test problems. Numerical results are included.
\end{abstract}
{\bf Keywords}: Lanczos Algorithm; Systems of Linear Equations; Formal Orthogonal Polynomials.

\bigskip
\noindent {\bf 2010 Mathematical Subject Classification: 65F10}

\section{Introduction}
The Lanczos algorithm, \cite{50:Lanczos,04:Broyden}, has been designed to find the eigenvalues of a matrix. However, it has found application in the area of Systems of Linear Equations (SLE's) where it is now well established. It is an iterative process which, in exact arithmetic, finds the exact solution in at most $n$ number of steps \cite{52:Lanczos}, where $n$ is the dimension of the problem. Several Lanczos-type algorithms have been designed and among them, the famous conjugate gradient algorithm of Hestenes and Stiefel \cite{52:Hest}, when the matrix is Hermitian and the bi-conjugate gradient algorithm of Fletcher \cite{76:Fletcher}, in the general case. In the last few decades, Lanczos-type algorithms have evolved and different variants have been derived, which can be found in
\cite{94:Baheux,98:Bjorck,93:Brezinski,94:Brezinski,99:Brezinski,00:Brezinski,02:Brezinski,92:Brezinski,00:Calvetti,97:Greenbaum,99:Guennouni,06:Meurant,79:Parlett,85:Parlett,87:Saad,87:Van,94:Ye,10:Farooq}.

Lanczos-type algorithms are commonly derived using Formal Orthogonal Polynomials (FOP's), \cite{93:Brezinski}. The connection  between the Lanczos algorithm, \cite{52:Lanczos} and orthogonal polynomials, \cite{39:Szego} has been studied extensively in
\cite{94:Baheux,93:Brezinski,00:Brezinski,02:Brezinski,92:Brezinski,80:Brezinski,91:Zaglia,91:Brezinski,83:Draux}.

\subsection{Notation}
The notation introduced by Baheux, in \cite{94:Baheux,95:Baheux}, for recurrence relations with three terms is adopted here. It puts recurrence relations involving FOP's $P_k(x)$ (the polynomials of degree at most $k$ with regard to the linear functional $c$) and/or FOP's $P^{(1)}_k(x)$ (the monic polynomials of degree at most $k$ with regard to linear functional $c^{(1)}$, \cite{91:Brezinski}) into two groups: $A_i$ and $B_j$. Although relations $A_i$, when they exist, rarely lead to Lanczos-type algorithms on their own (the exceptions being $A_4$, \cite{94:Baheux,95:Baheux}, and $A_{12}$, \cite{10:Farooq}), relations $B_j$ never lead to such algorithms for obvious reasons. It is the combination of recurrence relations $A_i$ and $B_j$, denoted as $A_i/B_j$, when both exist, that lead to Lanczos-type algorithms. In the following we will refer to algorithms by the relation(s) that lead to them. Hence, we will have, potentially, algorithms $A_i$ and algorithms $A_i/B_j$, for some $i=1,2,\dots $ and some $j=1,2,\dots $.

The paper is organized as follows. In the next section, the background to the Lanczos process is presented. Section 3 is on FOP's. Section 4 is on algorithm $A_{12}$, \cite{10:Farooq} and the estimation of the coefficients of the recurrence relations $A_{12}$ used to derive it. Section 5 is the estimation of the coefficients of recurrence relation $A_{12}$, \cite{10:Farooq}, used to derive the new algorithm of the same name $i.e$ $A_{12}(new)$. Section 6 describes the test problems and reports the numerical results. Section 7 is the conclusion and further work.

\section{The Lanczos Process}
Consider the system of linear equations,
\begin{equation}
\textit{A}\textbf{x}=\textbf{b},
\end{equation}
\noindent where $\textit{A}$ is $n\times n$ real matrix, $\textbf{b}$ and $\textbf{x}$ are vectors in $\textit{R}^{n}.$

 Choose $\textbf{x}_{0}$ and $\textbf{y}$, two arbitrary vectors in $\textit{R}^{n}$, such that $\textbf{y}\neq0$. Then, Lanczos process    \cite{52:Lanczos} consists in generating a sequence of vectors $\textbf{x}_{k}\in \textit{R}^n$, such that
\begin{equation}\textbf{x}_{k}-\textbf{x}_{0}\in \textit{F}_{k}(\textit{A}, \textbf{r}_{0})=\texttt{span}(\textbf{r}_{0}, \textit{A}\textbf{r}_{0},\ldots,\textit{A}^{k-1}\textbf{r}_{0}),\end{equation}
and
\begin{equation}\textbf{r}_{k}=\textbf{b}-\textit{A}\textbf{x}_{k}\bot\textit{G}_{k}(\textit{A}^{T}, \textbf{y})=\texttt{span}(\textbf{y}, \textit{A}^{T}\textbf{y},\ldots,\textit({A}^{T})^{k-1}\textbf{y}),\end{equation}
\noindent where $\textit{A}^{T}$ is the transpose of matrix $\textit{A}$.

Equation $(2)$ implies,
\begin{equation}
\textbf{x}_{k}-\textbf{x}_{0}=-\beta_{1}\textbf{r}_{0}-\beta_{2}\textit{A}\textbf{r}_{0}-\dots-\beta_{k}\textit{A}^{k-1}\textbf{r}_{0}.
\end{equation}
\noindent Multiplying both sides by $\textit{A}$ and adding and subtracting $\textbf{b}$ on the left hand side of (4) gives
\begin{equation}
\textbf{r}_{k}=\textbf{r}_{0}+\beta_{1}\textit{A}\textbf{r}_{0}+\beta_{2}\textit{A}^{2}\textbf{r}_{0}+\dots+\beta_{k}\textit{A}^{k}\textbf{r}_{0}.
\end{equation}
\noindent If we set
\[P_{k}(x)=1+\beta_{1}x+\dots+\beta_{k}x^{k},\]
\noindent then we can write from $(5)$
\begin{equation}
\textbf{r}_{k}=P_{k}(\textit{A})\textbf{r}_{0}.
\end{equation}
\noindent From $(3)$, since $(\it{A^T})^{i}\bf{y}$ and $\bf{r}_k$ are each in orthogonal subspaces, we can write,\\

$((\it{A^T})^{i}\bf{y},\bf{r}_{k}) =(\bf{y}, \it{A}^i\bf{r}_{k})=(\bf{y}, \it{A}^i\it{P}_{k}(\it{A})\bf{r}_{0})=0$, $\mbox{ for } i=0,\dots,k-1.$\\

\noindent Thus, the coefficients $\beta_{1}$,\dots,$\beta_{k}$ form a solution of system of linear equations,
\begin{equation}
\beta_{1}(\textbf{y}, \textit{A}^{i+1}\textbf{r}_{0})+\dots+\beta_{k}(\textbf{y}, \textit{A}^{i+k}\textbf{r}_{0})=-(\textbf{y}, \textit{A}^{i}\textbf{r}_{0}), \mbox{ for } i=0,\dots,k-1.\end{equation}
\noindent If the determinant of the above system is not zero then its solution exists and allows to obtain $\textbf{x}_{k}$ and $\textbf{r}_{k}$. Obviously, in practice, solving the above system directly for increasing values of $k$ is not feasible; $k$ is the order of the iterate in the solution process.  We shall see now how to solve this system for increasing values of $k$ recursively, that is, if polynomials $P_{k}$ can be computed recursively. Such computation is feasible as the polynomials $P_{k}$ form a family of FOP's and will now be explained. In exact arithmetic, $k$ should not exceed $n$, where $n$ is the dimension of the problem.

\section{Formal Orthogonal Polynomials}
Let $c$ be a linear functional on the space of complex polynomials defined by
\begin{center}
$c(x^i)=c_{i}$ \mbox{ for } $i=0,1,\dots$
\end{center}

\noindent where \begin{center} $c_{i}=({(\textit{A}^T)}^{i}\textbf{y},\textbf{r}_{k})=(\textbf{y}, \textit{A}^i\textbf{r}_{k})$ \mbox{ for } $i=0,1,\dots$\end{center}
\noindent Again, because of (3) above, an orthogonality condition can be written as,
\begin{equation}
c(x^iP_{k})=0 \mbox{ for } i=0,\dots,k-1.
\end{equation}
\noindent This condition shows that $P_{k}$ is the polynomial of degree at most $k$ which is a FOP with respect to the functional $c$, \cite{80:Brezinski}.

Given the expression of $P_k(x)$ above, $P_{k}(0)=1, \forall k$ is a normalization condition for these polynomials; $P_{k}$ exists and is unique if the following Hankel determinant
\[\textit{H}^{(1)}_{k}=\left\vert\begin{array}{cccc}
c_{1} & c_{2} & \cdots & c_{k}\\
c_{2} & c_{3} & \cdots & c_{k+1}\\
\vdots & \vdots &  & \vdots\\
c_{k} & c_{k+1} & \cdots & c_{2k-1}\\
\end{array}\right\vert\]\noindent is not zero, in which case we can write $P_{k}(x)$ as
\begin{equation}
P_{k}(x)=\frac{\left\vert\begin{array}{cccc}
1 & x & \cdots & x^k\\
c_{0} & c_{1} & \cdots & c_{k}\\
\vdots & \vdots &  & \vdots\\
c_{k-1} & c_{k} & \cdots & c_{2k-1}\\
\end{array}\right\vert}{\left\vert\begin{array}{cccc}
c_{1} & \cdots & c_{k}\\
\vdots &  & \vdots\\
c_{k} & \cdots & c_{2k-1}\\
\end{array}\right\vert},
\end{equation}
\noindent where the denominator of this polynomial is $\textit{H}^{(1)}_{k}$, the determinant of the system (7). We assume that $\forall$ $k$, $\textit{H}^{(1)}_{k}\neq0$ and therefore all the polynomials $P_{k}$ exist for all $k$. If for some $k$, $\textit{H}^{(1)}_{k}=0$, then $P_{k}$ does not exist and breakdown occurs in the algorithm, \cite{93:Brezinski,00:Brezinski,02:Brezinski,92:Brezinski,91:Zaglia}.

A Lanczos-type method consists in computing $P_{k}$ recursively, then $\textbf{r}_{k}$ and finally $\textbf{x}_{k}$ such that $\textbf{r}_{k}=\textbf{b}-\textit{A}\textbf{x}_{k}$, without inverting $A$. This gives the solution of the system $(1)$ in at most $n$ steps, where $n$ is the dimension of the SLE. For more details, see \cite{93:Brezinski,99:Brezinski}.

FOPs can be put together into recurrence relations. Such relations give rise to various procedures for the recursive computation of $P_{k}$ and hence we get different Lanczos-type algorithms for computing $\textbf{r}_{k}$ and, therefore, $\textbf{x}_{k}$. These algorithms have been  studied in \cite{94:Baheux,93:Brezinski,99:Brezinski,00:Brezinski,02:Brezinski,92:Brezinski,91:Zaglia,95:Baheux}. They differ by the recurrence relationships used to express the polynomials $P_{k}, k=2,3,...$.

\section{Recurrence relation $A_{12} based algorithm$}

Algorithms $A_5/B_{10}$, $A_8/B_{10}$ and $A_{12}$ are the most robust algorithms as found in \cite{94:Baheux,95:Baheux,10:Farooq,11:Salhi,12:Farooq}, on the same problems considered here. We, therefore compare our results with these algorithms. Since the algorithm we introduce here is also based on the recurrence relation $A_{12}$ \cite{10:Farooq,11:Farooq}, according to the notation of \cite{94:Baheux}, it is really a modification of algorithm $A_{12}$ that can be found in \cite{10:Farooq}. Indeed, $A_{12}$ is derived using the auxiliary polynomial $U_i(x)=x^i$, of exact degree $i$, while here we derive a
new algorithm $A_{12}$ but for $U_i(x)=P_i(x)$. For completeness, we recall algorithm $A_{12}$ here but leave out its derivation which can be found by the interested reader in \cite{10:Farooq}.

\subsection{Algorithm $A_{12}$}
Algorithm $A_{12}$  \cite{10:Farooq} can be described as follows.
\begin{algorithm}[H]
\caption{: Lanczos-type algorithm $A_{12}$}
\begin{algorithmic}
\STATE Choose $\textbf{x}_{0}$ and $\textbf{y}$ such that $\textbf{y}\neq0$,

\STATE Set $r_{0}=b-Ax_{0}$, $y_{0}=y$, $p=Ar_{0}$, $p_{1}=Ap$,
$c_{0}=(y, r_{0})$, \STATE  $c_{1}=(y, p)$, $c_{2}=(y, p_{1})$,
$c_{3}=(y, Ap_{1})$, $\delta=c_{1}c_{3}-c_{2}^2$, \STATE
$\alpha=\frac{c_{0}c_{3}-c_{1}c_{2}}{\delta}$,
$\beta=\frac{c_{0}c_{2}-c_{1}^2}{\delta}$, \STATE
$r_{1}=r_{0}-\frac{c_{0}}{c_{1}}p$,
$x_{1}=x_{0}+\frac{c_{0}}{c_{1}}r_{0}$, \STATE
$r_{2}=r_{0}-\alpha p+\beta p_{1}$, $x_{2}=x_{0}+\alpha
r_{0}-\beta p$, \STATE $y_{1}=A^{T}y_{0}$, $y_{2}=A^{T}y_{1}$,
$y_{3}=A^{T}y_{2}$.

\FOR {k = 3, 4,\dots,} \STATE $y_{k+1}=A^{T}y_{k}$,
$q_{1}=Ar_{k-1}$, $q_{2}=Aq_{1}$, $q_{3}=Ar_{k-2}$, \STATE
$a_{11}=(y_{k-2}, r_{k-2})$, $a_{13}=(y_{k-3}, r_{k-3})$,
$a_{21}=(y_{k-1}, r_{k-2})$, $a_{22}=a_{11}$,\STATE
$a_{23}=(y_{k-2}, r_{k-3})$, $a_{31}=(y_{k}, r_{k-2})$,
$a_{32}=a_{21}$, $a_{33}=(y_{k-1}, r_{k-3})$, \STATE
$s=(y_{k+1}, r_{k-2})$, $t=(y_{k}, r_{k-3})$,
$F_{k}=-\frac{a_{11}}{a_{13}}$,\STATE
$b_{1}=-a_{21}-a_{23}F_{k}$, $b_{2}=-a_{31}-a_{33}F_{k}$,
$b_{3}=-s-tF_{k}$,\STATE
$\Delta_{k}=a_{11}(a_{22}a_{33}-a_{32}a_{23})+a_{13}(a_{21}a_{32}-a_{31}a_{22})$,
\STATE
$B_{k}=\frac{b_{1}(a_{22}a_{33}-a_{32}a_{23})+a_{13}(b_{2}a_{32}-b_{3}a_{22})}{\Delta_{k}}$,
\STATE $G_{k}=\frac{b_{1}-a_{11}B_{k}}{a_{13}}$, \STATE
$C_{k}=\frac{b_{2}-a_{21}B_{k}-a_{23}G_{k}}{a_{22}}$, \STATE
$A_{k}=\frac{1}{C_{k}+G_{k}}$, \STATE
$r_{k}=A_{k}\{q_{2}+B_{k}q_{1}+C_{k}r_{k-2}+F_{k}q_{3}+G_{k}r_{k-3}\}$,
\STATE
$x_{k}=A_{k}\{C_{k}x_{k-2}+G_{k}x_{k-3}-(q_{1}+B_{k}r_{k-2}+F_{k}r_{k-3})\}$;
\STATE {\bf If} $||r_{k}|| \leq \epsilon$, then $x=x_{k}$, {\bf Stop}.
\ENDFOR
\end{algorithmic}
\end{algorithm}

\section{The new algorithm $A_{12}$ and its derivation}
As said above, in \cite{10:Farooq}, relation $A_{12}$ is derived using the auxiliary polynomial $U_i(x)=x^i$, of exact degree $i$. Here, we discuss the same relation, but for $U_i(x)=P_i(x)$. Coefficients are estimated for the new case. The Lanczos-type algorithm based on $A_{12}$ for the new choice of $U_i$, is called $A_{12}(new)$. This new algorithm is described below. Before deriving and discussing it, we recall the definition of an orthogonal polynomials sequence, \cite{10:Farooq}.\\

\textbf{Definition} 1. A sequence ${P_m}$ is called an orthogonal polynomial sequence, \cite{84:Chihara} with respect to the linear functional $c$ if, for all nonnegative integers $n$ and $m$,\\
(i) $P_m$ is a polynomial of degree $m$,\\
(ii) $c(x^nP_m) = 0$, for $m \neq n$,\\
(iii) $c(x^mP_m)$ $\neq 0.$

\subsection{Relation $A_{12}$ for the choice $U_i(x)=P_i(x)$.}

Consider the following recurrence relationship, \cite{10:Farooq},
\begin{equation}\label{equa}
P_k(x)=A_k\{(x^2+B_kx+C_k)P_{k-2}(x)+(D_kx^3+E_kx^2+F_kx+G_k)P_{k-3}(x)\}\\
\end{equation}
where $P_k$, $P_{k-2}$, and $P_{k-3}$ are polynomials of degree $k$, $k-2$, and $k-3$, respectively. $A_k$, $B_k$, $C_k$, $D_k$, $E_k$, $F_k$ and $G_k$ are the coefficients to be determined using the normality and the orthogonality conditions given in Section 3. Let, again, $c$ be a linear functional defined by $c(x^i)=c_i$. The orthogonality condition gives
\[c(U_iP_k)=0, i=0,1\cdots, k-1.\]

For $x=0$, and applying the normality condition, (\ref{equa}) becomes
\begin{equation}
1=A_k\{C_k+G_k\}.
\end{equation}
Now multiply (\ref{equa}) by $U_i$. Applying `c' on both sides and using the orthogonality condition, we get
\begin{eqnarray}\label{5}
c(x^2U_iP_{k-2})+B_kc(xU_iP_{k-2})+C_kc(U_iP_{k-2})+D_kc(x^3U_iP_{k-3})\nonumber\\
+E_kc(x^2U_iP_{k-3})+F_kc(xU_iP_{k-3})+G_kc(U_iP_{k-3})=0.
\end{eqnarray}
The orthogonality condition holds for $i=0,1,2,\cdots,k-7.$\\
For $i=k-6$, equation $(\ref{5})$ gives
\begin{center}
$D_kc(x^3U_{k-6}P_{k-3})=0$,\\
\end{center}
which implies that $D_k=0$, since $c(x^3U_{k-6}P_{k-3}) \neq 0$.\\
For $i=k-5$, $(\ref{5})$ becomes $E_kc(x^2U_{k-5}P_{k-3})=0$.\\
Since $c(x^2 U_{k-5}P_{k-3})\neq0$, $E_k=0$.\\
For $i=k-4$, we get \\

$c(x^2U_{k-4}P_{k-2})+F_kc(xU_{k-4}P_{k-3})=0$, which gives
\begin{equation}\label{FU}
F_k=-\frac{c(x^2U_{k-4}P_{k-2})}{c(xU_{k-4}P_{k-3})}.\end{equation}
For $i=k-3$, $i=k-2$ and $i=k-1$  equation $(\ref{5})$ can be respectively written as,
\begin{eqnarray}\label{11}
B_kc(xU_{k-3}P_{k-2})+G_kc(U_{k-3}P_{k-3})=-c(x^2U_{k-3}P_{k-2})-F_kc(xU_{k-3}P_{k-3}),
\end{eqnarray}
\begin{eqnarray}\label{12}
B_kc(xU_{k-2}P_{k-2})+C_kc(U_{k-2}P_{k-2})+G_kc(U_{k-2}P_{k-3})=\nonumber\\
-c(x^2U_{k-2}P_{k-2})-F_kc(xU_{k-2}P_{k-3}),
\end{eqnarray}
\begin{eqnarray}\label{13}
B_kc(xU_{k-1}P_{k-2})+C_kc(U_{k-1}P_{k-2})+G_kc(U_{k-1}P_{k-3})=\nonumber\\
-c(x^2U_{k-1}P_{k-2})-F_kc(xU_{k-1}P_{k-3}).
\end{eqnarray}
Now for simplicity let us denote the right sides of equations $(\ref{11})$, $(\ref{12})$ and $(\ref{13})$ by $b_1$ ,$b_2$ and $b_3$ respectively then we get the following system of equations,\\
\begin{equation}
B_kc(xU_{k-3}P_{k-2})+G_kc(U_{k-3}P_{k-3})=b_1,
\end{equation}
\begin{equation}
B_kc(xU_{k-2}P_{k-2})+C_kc(U_{k-2}P_{k-2})+G_kc(U_{k-2}P_{k-3})=b_2,
\end{equation}
\begin{equation}
B_kc(xU_{k-1}P_{k-2})+C_kc(U_{k-1}P_{k-2})+G_kc(U_{k-1}P_{k-3})=b_3.
\end{equation}

If $\Delta_k$ denotes the determinant of the coefficient matrix of the above system of equations then, \\
\begin{eqnarray}\label{delta}
\Delta_k=c(xU_{k-3}P_{k-2})\{c(U_{k-2}P_{k-2})c(U_{k-1}P_{k-3})-c(U_{k-2}P_{k-3})c(U_{k-1}P_{k-2})\}+\nonumber\\c(U_{k-3}P_{k-3})\{c(xU_{k-2}P_{k-2})c(U_{k-1}P_{k-2})-c(U_{k-2}P_{k-2})c(xU_{k-1}P_{k-2})\}.
\end{eqnarray}
If $\Delta_k\neq0$ then,\\
\begin{eqnarray*}
B_k=\frac{b_1}{\Delta_k}\{c(U_{k-2}P_{k-2})c(U_{k-1}P_{k-3})-c(U_{k-2}P_{k-3})c(U_{k-1}P_{k-2})\}\\
\\ +\frac{c(U_{k-3}P_{k-3})\{b_2c(U_{k-1}P_{k-2})-b_3c(U_{k-2}P_{k-2})\}}{\Delta_k},
\end{eqnarray*}
\begin{eqnarray*}
G_k=\frac{b_1-c(xU_{k-3}P_{k-2})B_k}{c(U_{k-3}P_{k-3})},
\end{eqnarray*}
\begin{eqnarray*}
C_k=\frac{b_2-c(xU_{k-2}P_{k-2})B_k-c(U_{k-2}P_{k-3})G_k}{c(U_{k-2}P_{k-2})}.
\end{eqnarray*}
With the above new estimated coefficients, the expression of polynomials $P_k(x)$ can be written as,
\begin{equation}\label{s}
P_k(x)=A_k\{(x^2+B_kx+C_k)P_{k-2}(x)+(F_kx+G_k)P_{k-3}(x)\}.
\end{equation}
Now, for $U_i(x)=P_k(x)$, and from equation (\ref{FU}), $F_k$  becomes\\
\begin{eqnarray*}F_k=-\frac{c(x^2P_{k-4}P_{k-2})}{c(xP_{k-4}P_{k-3})}.\end{eqnarray*}
Similarly, from equation (\ref{delta}) $\Delta_k$ becomes,\\
$\Delta_k=c(xP_{k-3}P_{k-2})\{c(P_{k-2}^2)c(P_{k-1}P_{k-3})-c(P_{k-2}P_{k-3})c(P_{k-1}P_{k-2})\}+\\c(P_{k-3}^2)\{c(xP_{k-2}^2)c(P_{k-1}P_{k-2})-c(P_{k-2}^2)c(xP_{k-1}P_{k-2})\}.$\\

Using (definition $ 1$), \cite{10:Farooq}, $\Delta_k$ simplifies to \\
\begin{eqnarray*}\Delta_k=-c(P_{k-3}^2)c(P_{k-2}^2)c(xP_{k-1}P_{k-2}).\end{eqnarray*}
Using again definition $1$ and $U_i(x)=P_k(x)$,  the rest of the coefficients can be determined as follows. Let\\
\[b_1=-c(x^2P_{k-3}P_{k-2})-F_kc(xP_{k-3}^2),\]
\[b_2=-c(x^2P_{k-2}^2)-F_kc(xP_{k-2}P_{k-3}),\]
\[b_3=-c(x^2P_{k-1}P_{k-2})-F_kc(xP_{k-1}P_{k-3}),\]
then\\
\begin{eqnarray*}B_k = \frac{1}{\Delta_k}\{b_1\{c(P_{k-2}P_{k-2})c(P_{k-1}P_{k-3})-c(P_{k-2}P_{k-3})c(P_{k-1}P_{k-2})\}\\+c(P_{k-3}P_{k-3})\{b_2c(P_{k-1}P_{k-2})-b_3c(P_{k-2}P_{k-2})\}\},\end{eqnarray*}
or, \begin{eqnarray*}B_k =  -\frac{b_3c(P_{k-3}^2)c(P_{k-2}^2)}{\Delta_k} = -\frac{b_3}{c(xP_{k-1}P_{k-2})},\end{eqnarray*}
\begin{eqnarray*}G_k = \frac{b_1-c(xP_{k-3}P_{k-2})B_k}{c(P_{k-3}^2)},\end{eqnarray*}
\begin{eqnarray*}C_k=\frac{b_2-c(xP_{k-2}P_{k-2})B_k-c(P_{k-2}P_{k-3})G_k}{c(P_{k-2}P_{k-2})}= \frac{b_2-c(xP_{k-2}^2)B_k}{c(P_{k-2}^2)},\end{eqnarray*} and
\begin{eqnarray*}A_k=\frac{1}{C_k+G_k}.\end{eqnarray*}
 As in \cite{10:Farooq}, we can write,
\begin{equation}
 \textbf{r}_k=A_k\{\textit{A}^2\textbf{r}_{k-2}+B_k \textit{A}\textbf{r}_{k-2}+C_k \textbf{r}_{k-2}+F_k \textit{A}\textbf{r}_{k-3}+G_k \textbf{r}_{k-3}\},
 \end{equation}
 \begin{equation}
 \textbf{x}_k=A_k\{C_k\textbf{x}_{k-2}+G_k \textbf{x}_{k-3}-(\textit{A}\textbf{r}_{k-2}+B_k \textbf{r}_{k-2}+F_k\textbf{r}_{k-3})\}.
\end{equation}
As we know from \cite{94:Baheux,95:Baheux},\\ 
\begin{eqnarray}\label{C}\begin{cases} $ setting $  U_k(x)=P_k(x) $ and $ \textbf{z}_k=P_k(A^T)\textbf{y}, $ we get$ \\
c(U_kP_k)=(y,U_k(A)P_(A)\textbf{r}_0)=(U_k(A^T)\textbf{y},P_k(A)\textbf{r}_0)=(\textbf{z}_k,\textbf{r}_k).\end{cases}
 \end{eqnarray}
  
  So, from relation (\ref{s}), after replacing $\textbf{x}$ by $A^T$, multiplying by $\textbf{y}$ on both sides and using (\ref{C})  we can write,
 \begin{equation} \textbf{z}_k=A_k\{(\textit{A}^T)^2\textbf{z}_{k-2}+B_k \textit{A}^T\textbf{z}_{k-2}+C_k \textbf{z}_{k-2}+F_k \textit{A}^T\textbf{z}_{k-3}+G_k \textbf{z}_{k-3}\}.
 \end{equation}
 Similarly using (\ref{C}) all coefficients become,\\
$F_k=-\frac{c(x^2P_{k-4}P_{k-2})}{c(xP_{k-4}P_{k-3})}$=$-\frac{(\textit{A}^T \textbf{z}_{k-2},\textit{A}\textbf{r}_{k-4})}{(\textbf{z}_{k-3},\textit{A}\textbf{r}_{k-4})}$,\\ 
 $\Delta_k=-c(P_{k-3}^2) c(P_{k-2}^2)c(xP_{k-1}P_{k-2})$=$-(\textbf{z}_{k-3},\textbf{r}_{k-3})(\textbf{z}_{k-2},\textbf{r}_{k-2})(\textbf{z}_{k-1},\textit{A}\textbf{r}_{k-2})$.\\
  $b_1=-(A^T \textbf{z}_{k-3},A\textbf{r}_{k-2})- F_k(\textbf{z}_{k-3},A \textbf{r}_{k-3}),$\\
  $b_2=-(A^T \textbf{z}_{k-2},A\textbf{r}_{k-2})- F_k(\textbf{z}_{k-2},A \textbf{r}_{k-3}),$\\
  $b_3=-(A^T\textbf{z}_{k-1},A\textbf{r}_{k-2})- F_k(\textbf{z}_{k-1},A \textbf{r}_{k-3}),$\\  $B_k$ = $\frac{b_3}{c(xP_{k-1}P_{k-2})}$ = $\frac{b_3}{(\textbf{z}_{k-1},A\textbf{r}_{k-2})}$,\\\\
  $G_k$=$\frac{b_1-c(xP_{k-3}P_{k-2})B_k}{c(P_{k-3}^2)}$=$\frac{b_1-(\textbf{z}_{k-3},A\textbf{r}_{k-2})B_k}{(\textbf{z}_{k-3},\textbf{r}_{k-3})}$,\\\\ $C_k$=$\frac{b_2-c(xP_{k-2}^2)B_k}{c(P_{k-2}^2)}$=$\frac{b_2-(\textbf{z}_{k-2},A\textbf{r}_{k-2})B_k}{(\textbf{z}_{k-2},\textbf{r}_{k-2})}$,\\\\
  $A_k$ = $\frac{1}{C_k+G_k}$.\\
  All previous formulae are valid for $k\geq 4$. So we need $\textbf{r}_1$, $\textbf{r}_2$, $\textbf{r}_3$ and $\textbf{z}_1$, $\textbf{z}_2$, $\textbf{z}_3$ to calculate $\textbf{r}_k$ and $\textbf{z}_k$ recursively. $\textbf{r}_1$ , $\textbf{r}_2$ and $\textbf{z}_1$, $\textbf{z}_2$ are found differently in \cite{10:Farooq}, while $\textbf{r}_3$ and $\textbf{z}_3$ can be determined in a similar way giving,\\
$\textbf{r}_3=\textbf{r}_0-\frac{\acute{\alpha}}{\Delta} \textbf{p}+\frac{\acute{\beta}}{\Delta}\textbf{p}_1-\frac{\acute{\gamma}}{\Delta} \textbf{p}_2$,\\
   $\textbf{z}_3=\textbf{z}_0-\frac{\acute{\alpha}}{\Delta} \textbf{y}_1+\frac{\acute{\beta}}{\Delta}\textbf{y}_2-\frac{\acute{\gamma}}{\Delta} \textbf{y}_3$.\\
   Using $\textbf{r}_k=\textbf{b}-\emph{A}\textbf{x}_k$, we get from $\textbf{r}_3$,\\
   $\textbf{x}_3=\textbf{x}_0+\frac{\acute{\alpha}}{\Delta} \textbf{r}_0
-\frac{\acute{\beta}}{\Delta}\textbf{p}+\frac{\acute{\gamma}}{\Delta} \textbf{p}_1$,
   where\\
   $\Delta=c_1(c_3c_5-c_4^{2})-c_2(c_2c_5-c_3c_4)+c_3(c_2c_4-c_3^{2})$,\\
   $\acute{\alpha}=c_0(c_3c_5-c_4^{2})-c_2(c_1c_5-c_2c_4)+c_3(c_1c_4-c_3c_2)$,\\
   $\acute{\beta}=c_0(c_2c_5-c_4c_3)-c_1(c_1c_5-c_2c_4)+c_3(c_1c_3-c_2^{2})$,\\
   $\acute{\gamma}=c_0(c_2c_4-c_3^{2})-c_1(c_1c_4-c_2c_3)+c_2(c_1c_3-c_2^{2})$.\\
Note that parameters {\bf{p},\bf{p}$_1$,\bf{p}$_2$,\bf{y}$_1$,\bf{y}$_2$,\bf{y}$_3$}, $\Delta, \alpha, \beta,$ and $\gamma$ are temporary and defined in the algorithm below.

\subsection{Algorithm $A_{12}(new)$}
    We can now describe the new variant of algorithm $A_{12}(new)$ as follows.
 \begin{algorithm}[H]
\caption{: Lanczos-type Algorithm $A_{12}(new)$.}
\begin{algorithmic}
\STATE  Choose $\textbf{x}_0$ and $\textbf{y}$ such that  $\textbf{y}\neq0$.\\
  Set $\textbf{r}_0=\textbf{b}-A\textbf{x}_0$, $\textbf{z}_0=\textbf{y}$,  $\textbf{p}=A \textbf{r}_0$, $\textbf{p}_1=A \textbf{p}$, $\textbf{p}_2=A \textbf{p}_1$, $\textbf{p}_3=A\textbf{p}_2$, $\textbf{p}_4=A \textbf{p}_3$,\\
   $c_0=(\textbf{y},\textbf{r}_0)$, $c_1=(\textbf{y},\textbf{p})$,  $c_2=(\textbf{y},\textbf{p}_1)$, $c_3=(\textbf{y},\textbf{p}_2)$, $ c_4=(\textbf{y},\textbf{p}_3)$, $c_5=(\textbf{y},\textbf{p}_4)$,\\
   $\delta=c_1c_3-c_2^2$, $\alpha=\frac{c_0c_3-c_1c_2}{\delta}$, $\beta=\frac{c_0c_2-c_1^2}{\delta}$,\\
   $\textbf{r}_1=\textbf{r}_0-(\frac{c_0}{c_1})\textbf{p}$, $\textbf{x}_1=\textbf{x}_0+(\frac{c_0}{c_1})\textbf{r}_0$,\\
   $\textbf{r}_2=\textbf{r}_0-\alpha \textbf{p}+\beta \textbf{p}_1$, $\textbf{x}_2=\textbf{x}_0+\alpha  \textbf{r}_0-\beta \textbf{p}$,\\
   $\textbf{y}_1=A^T\textbf{y}$, $\textbf{y}_2=A^T \textbf{y}_1$, $\textbf{y}_3=A^T \textbf{y}_2$,\\
   $\textbf{z}_1=\textbf{z}_0-(\frac{c_0}{c_1})\textbf{y}_1$, $\textbf{z}_2=\textbf{z}_0-\alpha \textbf{y}_1+\beta  \textbf{y}_2$,\\
   $\Delta=c_1(c_3c_5-c_4^{2})-c_2(c_2c_5-c_3c_4)+c_3(c_2c_4-c_3^{2})$,\\
   $\acute{\alpha}=c_0(c_3c_5-c_4^{2})-c_2(c_1c_5-c_2c_4)+c_3(c_1c_4-c_3c_2)$,\\
   $\acute{\beta}=c_0(c_2c_5-c_4c_3)-c_1(c_1c_5-c_2c_4)+c_3(c_1c_3-c_2^{2})$,\\
   $\acute{\gamma}=c_0(c_2c_4-c_3^{2})-c_1(c_1c_4-c_2c_3)+c_2(c_1c_3-c_2^{2})$,\\
   $\textbf{r}_3=\textbf{r}_0-\frac{\acute{\alpha}}{\Delta} \textbf{p}+\frac{\acute{\beta}}{\Delta}\textbf{p}_1-\frac{\acute{\gamma}}{\Delta} \textbf{p}_2$,\\
   $\textbf{z}_3=\textbf{z}_0-\frac{\acute{\alpha}}{\Delta} \textbf{y}_1+\frac{\acute{\beta}}{\Delta}\textbf{y}_2-\frac{\acute{\gamma}}{\Delta} \textbf{y}_3$,\\
   $\textbf{x}_3=\textbf{x}_0+\frac{\acute{\alpha}}{\Delta} \textbf{r}_0
-\frac{\acute{\beta}}{\Delta}\textbf{p}+\frac{\acute{\gamma}}{\Delta} \textbf{p}_1$.\\
\FOR {k = 4,5\dots,}
\STATE   $ q_1=A \textbf{r}_{k-2}$, $q_2=A \textbf{q}_1$, $q_3=A \textbf{r}_{k-3}$,
\STATE   $ \textbf{s}_1=A^T \textbf{z}_{k-2}$, $\textbf{s}_2=A^T \textbf{s}_1$, $\textbf{s}_3=A^T \textbf{z}_{k-3}$,
\STATE   $\Delta_k=-(\textbf{z}_{k-3},\textbf{r}_{k-3})(\textbf{z}_{k-2},\textbf{r}_{k-2})(\textbf{z}_{k-1},A\textbf{r}_{k-2})$,
\STATE   $F_k=-\frac{(A^T \textbf{z}_{k-2},A\textbf{r}_{k-4})}{(\textbf{z}_{k-3},A\textbf{r}_{k-4})}$,
\STATE   $b_{1}=-(A^T \textbf{z}_{k-3},A\textbf{r}_{k-2})- F_k(\textbf{z}_{k-3},A \textbf{r}_{k-3})$,
\STATE   $b_2=-(A^T \textbf{z}_{k-2},A\textbf{r}_{k-2})- F_k(\textbf{z}_{k-2},A \textbf{r}_{k-3}),$
\STATE   $b_3=-(A^T \textbf{z}_{k-1},A\textbf{r}_{k-2})- F_k(\textbf{z}_{k-1},A \textbf{r}_{k-3})$,
\STATE   $B_{k}=\frac{b_3}{(\textbf{z}_{k-1},A\textbf{r}_{k-2})}$ ,
\STATE   $G_k=\frac{b_1-(\textbf{z}_{k-3},A\textbf{r}_{k-2})B_k}{(\textbf{z}_{k-3},\textbf{r}_{k-3})}$,
\STATE   $C_k=\frac{b_2-(\textbf{z}_{k-2},A\textbf{r}_{k-2})B_k}{(\textbf{z}_{k-2},\textbf{r}_{k-2})}$,
\STATE   $A_k=\frac{1}{C_k+G_k}$,
\STATE   $\textbf{r}_k=A_k\{\textbf{q}_2+B_k \textbf{q}_1+C_k \textbf{r}_{k-2}+F_k \textbf{q}_{3}+G_k \textbf{r}_{k-3}\}$,
\STATE   $x_k=A_k\{C_k\textbf{x}_{k-2}+G_k \textbf{x}_{k-3}-(\textbf{q}_1+B_k \textbf{r}_{k-2}+F_k\textbf{r}_{k-3})\}$,
\STATE   $\textbf{z}_k=A_k\{\textbf{s}_2+B_k \textbf{s}_1+C_k \textbf{z}+F_k \textbf{s}_{3}+G_k \textbf{z}_{k-3}\}$. 
\STATE {\bf If} $||r_{k}|| \leq \epsilon$, then $x=x_{k}$, {\bf Stop}.
\ENDFOR
\end{algorithmic}
\end{algorithm}

\section{Numerical Tests}

$A_{12}(new)$ has been
tested against $A_{12}$, $A_5/B_{10}$ and $A_8/B_{10}$, the
best Lanczos-type algorithms according to \cite{94:Baheux,10:Farooq,95:Baheux}. The test problems arise in the 5-point discretisation of the operator $\frac{-\partial^{2}}{\partial x^2}-\frac{\partial^{2}}{\partial y^2}+\gamma\frac{\partial}{\partial x}$ on a rectangular region \cite{95:Baheux}. Comparative results on instances of the following problem ranging from dimension $10$ to $100$ for parameter $\delta$ taking values $0.0$  and for the tolerance $eps=1.0e-05$, are recorded in Table 1.

\[A=\left(\begin{array}{ccccccc}
B & -I & \cdots & \cdots  & 0\\
-I & B & -I  &  & \vdots\\
\vdots & \ddots & \ddots & \ddots & \vdots\\
\vdots &  & -I & B & -I\\
0 & \cdots & \cdots & -I & B\\
\end{array}
\right),\]

\noindent with
\[B=\left(\begin{array}{ccccccc}
4 & \alpha & \cdots & \cdots & 0\\
\beta & 4 & \alpha &  & \vdots\\
\vdots & \ddots & \ddots & \ddots & \vdots\\
\vdots & & \beta & 4 & \alpha\\
0 & \cdots &  & \beta & 4\\
\end{array}\right),\]

\noindent and $\alpha=-1+\delta$, $\beta=-1-\delta$. The right
hand side $\textbf{b}$ is taken to be
$\textbf{b}=\textit{A}\textbf{X}$, where $\textbf{X}=(1, 1,
\dots, 1)^{T}$, is the solution of the system. The dimension of
$\textit{B}$ is $10$.

\begin{table}[h]
\caption{Experimental results for problems when $\delta=0$}
\begin{center} 
\scalebox{0.7}{
\begin{tabular}{|c|c|c|c|c|c|c|c|c|}
\hline
  & \multicolumn{2}{c|}{$A_{5}/B_{10}$} &
\multicolumn{2}{c|}{$A_{8}/B_{10}$} & \multicolumn{2}{c|}{$A_{12}$} & \multicolumn{2}{c|}{$A_{12} new$}\\\cline{2-9}
$n$  & \multicolumn{1}{c}{$||r_{k}||$} & \multicolumn{1}{|c|}{t(sec)} &
\multicolumn{1}{c}{$||r_{k}||$} & \multicolumn{1}{|c|}{t(sec)} & \multicolumn{1}{c}{$||r_{k}||$} &
\multicolumn{1}{|c|}{t(sec)} & \multicolumn{1}{c}{$||r_{k}||$} & \multicolumn{1}{|c|}{t(sec)} \\\hline
$10$  & $2.0866e-013$ & 0.002628 & $3.5775e-013$ & 0.008440 & $1.0252e-013$ & 0.042433 & $2.7146e-015$&0.018560\\
$20$ & $ 2.5278e-014$ &0.002619 & $1.6765e-013$ &0.008624 & $1.8456e-013$ & $0.042880$& $2.4416e-015$& 0.017902\\
$30$  & $2.4011e-009$ &0.003139 & $6.9352e-009$ &0.009134 & $1.6272e-007$ & 0.043438& $ 2.0829e-010 $ &  0.019099\\
$40$  & $1.5539e-009$ &0.003344 & $1.5680e-009$ & 0.009113 & $2.0343e-010$ & 0.043924& $ 2.7946e-011$ &0.019164\\
$50$  & $1.8730e-006$ & 0.003810 & $1.4671e-006$ &0.009634 & $4.7570e-005$ & 0.044461 &$ 1.2734e-006$ &0.020314\\
$60$  & $5.9083e-006$ & 0.003747 & $6.6800e-006$ & 0.009599 & $2.8615e-005$ &0.044002& $ 2.3608e-006$ &0.020202\\
$70$ & $ 9.3260e-006$ &0.004658 & $ 4.6961e-006$ &0.010246 & $8.5638e-005$ & 0.044369& $5.3790e-007$ & 0.020988\\
$80$  & $ 4.5674e-006$ &0.005496 & $ 4.6144e-006$ & 0.011470 & $6.8618e-005$ &0.046109& $ 3.5468e-006$ &0.022625\\
$90$  & $NaN$ &  & $ NaN$ &  & $7.2121e-005$ &0.047276& $ 4.3695e-006$ &0.021556\\
$100$  & $ 9.0038e-006$ & 0.004284  & $8.4881e-007$ & 0.010383 & $3.1098e-005$ & 0.044758& $2.0040e-008$ &0.020606\\
\hline
\end{tabular}}
\end{center}
\label{turns}
\end{table}

Table 1 records the computational results obtained with algorithms $A_{12}$(new), $A_{12}$, $A_5/B_{10}$ and $A_8/B_{10}$. Clearly, $A_{12}$(new) is an improvement on
$A_{12}$ on both robustness/stability and efficiency accounts. Compared to the well established $A_5/B_{10}$ and $A_8/B_{10}$, it is definitely more robust/stable; indeed, all problems have been solved to the required accuracy by $A_{12}$(new), and the other two algorithms failed to do so in one case as evidenced by the "NaN" outputs which point to breakdown or lack of robustness and stability, on the problem of dimension n=90. On efficiency, however, as expected, algorithms $A_5/B_{10}$ and $A_8/B_{10}$ are faster since they rely on recurrence relations involving lower order FOP's requiring few coefficients to estimate; unlike $A_{12}$ and $A_{12}$(new).
\section{Conclusion }

In this paper we have shown that, if the recurrence relation $A_{12}$ \cite{10:Farooq}, is determined for the choice of $U_i(x)=P_i(x)$, other than $x^i$ which is discussed in \cite{10:Farooq}, then a more robust algorithm $A_{12}(new)$ can be derived. The numerical performance of this algorithm compares well to that of three existing Lanczos-type algorithms, which were found to be the best among a number of such algorithms, \cite{94:Baheux,10:Farooq,95:Baheux}, on the same set of problems as considered here. Another achievement of $A_{12}(new)$ is that it can solve the above problem when its dimension is up to 500, while the rest of algorithms give results for problems with dimensions less or equal to 100.



\begin{thebibliography}{99}
\bibitem{50:Lanczos}C. Lanczos. An Iteration Method for the Solution of the Eigenvalue Problem of Linear Differential and Integeral Operators. {\em Journal of Research of the National Bureau of Standards}, {\bf 45}, {(1950)}, 255--282.

\bibitem{04:Broyden} C. G. Broyden and M. T. Vespucci. {\em Krylov Solvers For Linear Algebraic Systems}, Elsevier, Amsterdam, The Netherlands, 2004.

\bibitem{52:Lanczos} C. Lanczos. Solution of systems of linear equations by minimized iteration. {\em Journal of the National Bureau of Standards}, {\bf 49}, {(1952)}, 33--53.

\bibitem{52:Hest} M.R. Hestenes and E. Stiefel. Mehtods of conjugate gradients for solving linear systems. {\em Journal of the National Bureau of Standards}, {\bf 49}, {(1952)}, 409--436.

\bibitem{76:Fletcher} R. Fletcher, Conjugate gradient methods for indefinite systems. in: G.A. Watson, (Ed), {\em Numerical Analysis, Dundee 1975}, Lecture Notes in Mathematics,, volume 506. Springer, Berlin, 1976.

\bibitem{94:Baheux} C. Baheux, {\em Algorithmes d'implementation de la m\'{e}thode de Lanczos}, PhD thesis, University of Lille 1, France, 1994.


\bibitem{98:Bjorck} A. Bj\^{o}rck, T. Elfving, and Z. Strakos, Stability of Conjugate Gradient and Lanczos Methods for Linear Least Squares Problems, {\em SIAM Journal of Matrix Analysis and Application}, {\bf 19}, {(1998)}, 720--736.

\bibitem{93:Brezinski} C. Brezinski and H. Sadok, Lanczos-type algorithms for solving systems of linear equations, {\em Applied Numerical Mathematics}, {\bf 11}, {(1993)}, 443--473.

\bibitem{94:Brezinski} C. Brezinski and M. R. Zaglia, Hybird procedures for solving linear systems, {\em Numerische Mathematik}, {\bf 67}, {1994}, 1--19.

\bibitem{99:Brezinski} C. Brezinski, M. R. Zaglia, and H. Sadok. New look-ahead Lanczos-type algorithms for linear systems, {\em Numerische Mathematik}, {\bf 83}, {(1999)}, 53--85.

\bibitem{00:Brezinski} C. Brezinski, M. R. Zaglia, and H. Sadok. The matrix and polynomial approaches to Lanczos-type algorithms, {\em Journal of Computational and Applied Mathematics}, {\bf 123}, {(2000)}, 241--260.

\bibitem{02:Brezinski} C. Brezinski, M. R. Zaglia, and H. Sadok. A review of formal orthogonality in Lanczos-based methods. {\em Journal of Computational and Applied Mathematics}, {\bf 140}, {(2002)}, 81--98.

\bibitem{92:Brezinski} C. Brezinski, M. R. Zaglia, and H. Sadok. A Breakdown-free Lanczos type algorithm for solving linear systems, {\em Numerische Mathematik}, {\bf 63}, {(1992)}, 29--38.

\bibitem{00:Calvetti} D. Calvetti, L. Reichel, F. Sgallari, and G. Spaletta, A Regularizing Lanczos iteration method for underdetermined linear systems, {\em Jouranl of Computational and Applied Mathematics}, {\bf 115}, {(2000)} 101--120.

\bibitem{97:Greenbaum} A. Greenbaum. {\em Iterative Methods for Solving Linear System}, Societ for Industrial and Applied Mathematics, Philadelphia, 1997.

\bibitem{99:Guennouni} A. El Guennouni. A unified approach to some strategies for the treatment of breakdown in Lanczos-type algorithms. {\em Applicationes Mathematicae}, {\bf 26}, {(1999)}, 477--488.

\bibitem{06:Meurant} G. Meurant. {\em The Lanczos and conjugate gradient algorithms, From Theory to Finite Precision Computations}. SIAM, Philadelphia, 2006.

\bibitem{79:Parlett} B. N. Parlett and D. S. Scott. The Lanczos Algorithm With Selective Orthogonaliztion. {\em Mathematics of Computation}, {\bf 33}, {(1979)}, 217--238.
\bibitem{85:Parlett} B. N. Parlett, D. R. Taylor, and Z. A. Liu. A Look-Ahead Lanczos Algorithm for Unsymmetric Matrices. {\em Mathematics of Computation}, {\bf 44}, {(1985)}, 105--124.

\bibitem{87:Saad} Y. Saad. On the Lanczos method for solving linear system with several right-hand sides. {\em Mathematics of Computation}, {\bf 48}, {(1987)}, 651--662.
\bibitem{87:Van} H. A. Van der Vorst, An iterative solution method for solving f(\textit{A})\textbf{x}=\textbf{b}, using Krylov subspace information obtained for the symmetric positive definite matrix \textit{A}, {\em Journal of Computational and Applied Mathematics}, {\bf 18(2)}, {(1987)}, 249--263.

\bibitem{94:Ye} Q. Ye, A Breakdown-Free Variation of the Nonsymmetric Lanczos Algorithms, {\em Mathematics of Computation}, {\bf 62}, {(1994)}, 179--207.

\bibitem{10:Farooq} M. Farooq. {\em New Lanczos-type Algorithms and their Implementation}, PhD thesis, University of Essex, UK, 2011.

\bibitem{11:Farooq} M. Farooq and A. Salhi. New Recurrence Relationships Between Orthogonal Polynomials Which Lead to New Lanczos-type Algorithms, {\em Journal of Prime Research in Mathematics}, {\bf 8}, {(2012)}, 61-75.
    
\bibitem{11:Salhi} M. Farooq and A. Salhi. A {S}witching {A}pproach to {A}void {B}reakdown in {L}anczos-type {A}lgorithms,
{\em to appear in Applied Mathematics and Information Sciences in 2013.}
        
\bibitem{12:Farooq} M. Farooq and A. Salhi. A {P}re-emptive {R}estarting {A}pproach to {B}eating the {I}nherent {I}nstability of {L}anczos-type {A}lgorithms,
{\em to appear in Iranian Journal of Science and Technology, Transaction-A, Science, in 2013.}

\bibitem{39:Szego} G. Szeg$\ddot{o}$. Orthogonal Polynomials, {\em American Mathematical Society}, Providence, Rhode Island, 1939.

\bibitem{80:Brezinski} C. Brezinski, {\em Pad\'{e}-Type Approximation and General Orthogonal Polynomials}, Internat. Ser. Nuner. Math. 50. Birkh$\ddot{a}$user, Basel, 1980.

\bibitem{91:Zaglia} C. Brezinski and M. R. Zaglia, A new presentation of orthogonal polynomials with applications to their computation, {\em Numerical Algorithms}, {\bf 1}, {(1991)}, 207--222.

\bibitem{91:Brezinski} C. Brezinski, M. R. Zaglia, and H. Sadok. Avoiding breakdown and nearbreakdown in Lanczos type algorithms, {\em Numerical Algorithms}, {\bf 1}, {(1991)}, 261--284.

\bibitem{83:Draux} A. Draux, {\em Polyn\^{o}mes Orthogonaux Formels}, Application, LNM 974. Springer-Verlag, Berlin, 1983.

\bibitem{95:Baheux} C. Baheux, New Implementations of Lanczos Method. {\em Journal of Computational and Applied Mathematics}, {\bf 57}, {(1995)}, 3--15.

\bibitem{84:Chihara} T. S. Chihara, \textit{An Introduction to Orthogonal Polynomials}, Gordon and Breach, New York, London, Paris, 1984.

\end{thebibliography}
\end{document}